\documentclass[11pt,a4paper,psamsfonts,twoside]{amsart}
\usepackage{latexsym, amssymb,amscd,amsmath,amsthm}
\usepackage{graphicx} 

\chardef\bslash=`\\ 

\newtheorem{thm}{Theorem}[section]
\newtheorem{cor}[thm]{Corollary}
\newtheorem{lemma}[thm]{Lemma}
\newtheorem{prop}[thm]{Proposition}
\newtheorem{conj}[thm]{Conjecture}

\theoremstyle{remark}

\theoremstyle{definition}

\newtheorem*{ackn}{Acknowledgement}

\numberwithin{equation}{section}

\def\ubrace#1#2{\underbrace{#1}\limits_{\displaystyle{#2}}}
\def\Cal#1{{\mathcal#1}}
\def\<{\langle}\def\>{\rangle}
\def\what{\widehat}
\def\Z{{\mathbb Z}}\def\N{{\mathbb N}} 
\def\R{{\mathbb R}} 
\def\D{{\mathbb D}}
\newcommand{\xe}{\ensuremath{{x_\epsilon}}}
\newcommand{\Ce}{\ensuremath{{C^n_{\epsilon}}}}

\def\min{\textsl{min}}
\def\Mod{\textsl{Mod}}

\def\Fix{\textsl{Fix}}

\newcommand{\inv}{^{-1}}

\def\ga{\gamma}                 \def\Ga{\Gamma}

\begin{document}

\title[Relative hyperbolicity and Artin groups]
{Relative hyperbolicity and Artin groups}

\date{November 2005}
\author[R.~Charney]{Ruth Charney}
	\address{Department of Mathematics\\
	Brandeis University\\
	Waltham, MA 02454\\
	USA}
	\email{charney@brandeis.edu}
	\thanks{Charney was partially
      supported by NSF grant DMS-0405623.}
\author[J.~Crisp]{John Crisp}
	\address{I.M.B.(UMR 5584 du CNRS)\\ 
              Universit\'e de Bourgogne\\
		     B.P. 47 870\\
		21078 Dijon, France}
	\email{jcrisp@u-bourgogne.fr}
\thanks{ \hfill Typeset by \AmS-\LaTeX}

\keywords{relative hyperbolicity, Artin group, Deligne complex}

\subjclass{20F36}

\begin{abstract}
This paper considers the question of relative hyperbolicity of an Artin group 
with regard to the geometry of its associated Deligne complex. We prove 
that an Artin group is weakly hyperbolic relative to its finite (or spherical) 
type parabolic subgroups if and only if its Deligne complex is a Gromov hyperbolic space. 
For a 2-dimensional Artin group the Deligne complex is Gromov
hyperbolic precisely when the corresponding Davis complex is Gromov hyperbolic, 
that is, precisely when the underlying Coxeter group is a hyperbolic group.   
For Artin groups of FC type we give a sufficient condition for hyperbolicity of 
the Deligne complex which applies to a large class of these groups for which the 
underlying Coxeter group is hyperbolic.
\end{abstract}

\maketitle


\section{Introduction}\label{Intro}


Let $G$ denote a finitely generated group, and $\Cal H=\{ H_1, H_2,\dots, H_n\}$
a finite family of subgroups of $G$. Let $\Ga_S$ denote the Cayley graph of $G$ 
with respect to a finite generating set $S$. We denote by $\Ga_{S,\Cal H}$ the 
\emph{coned-off Cayley graph} with respect to $\Cal H$. Namely, $\Ga_{S,\Cal H}$
is the graph obtained from $\Ga_S$ by introducing a vertex $V_{gH}$ for each left coset $gH$,
with $g\in G$ and $H\in\Cal H$, and attaching $V_{gH}$ by an edge of length $\frac{1}{2}$
to each vertex of $\Ga_S$ labelled by an element of the coset $gH$. The isometric 
(left) action of $G$ on the Cayley graph $\Ga_S$ clearly extends to an isometric action of 
$G$ on $\Ga_{S,\Cal H}$ (by setting $g(V_{g'H})=V_{gg'H}$ for all $g,g'\in G$ and $H\in\Cal H$).
The nontrivial vertex stabilizers of this action are all conjugate to subgroups of $\Cal H$.

The group $G$ is said to be \emph{weakly hyperbolic relative to $\Cal H$}
if the coned-off Cayley graph $\Ga_{S,\Cal H}$ is a Gromov hyperbolic space.
This definition can be shown to be independent 
of the choice of finite generating set $S$ (see \cite{Farb}).
(We also refer the reader to \cite{Gro1,Gro2} or other standard 
references, such as \cite{BH}, for the definition of a Gromov hyperbolic space). 

The notion of weak relative hyperbolicity just described was first studied by B. Farb 
in \cite{Farb}.  We remark that the more conventional notion of relative hyperbolicity introduced
by M. Gromov \cite{Gro1,Gro2}, and explored by B. Bowditch \cite{Bow}, A. Yaman \cite{Yam} 
and others is rather more specific. 
It is equivalent to weak relative hyperbolicity as defined above
together with the additional condition of \emph{bounded coset penetration} (BCP) 
introduced by Farb in \cite{Farb}. The stronger form of relative hyperbolicity \`a la Gromov
models more precisely certain classical situations such as that of a geometrically 
finite Kleinian group
(which is hyperbolic relative to its parabolic subgroups in the stronger sense).

The difference between these two notions of relative hyperbolicity is perfectly and
rather starkly illustrated in the case of Artin groups 
(defined in Section \ref{Sect:Artin} below).
In  \cite{KS}, I. Kapovich and P. Schupp showed that  Artin groups 
with all relator indices at least $7$ are weakly hyperbolic relative 
to their non-free rank 2 parabolic subgroups (in this case the finite type
parabolics). On the other hand, as observed in \cite{KS}, the presence of numerous mutually
intersecting free abelian subgroups tends to rule out altogether the possibility of
relative hyperbolicity in the stronger sense for most Artin groups. 
More precisely, it is shown in \cite{BDM} and also follows easily
from Lemma 4 of \cite{AAS} that if an Artin group $G$ is strongly relatively hyperbolic 
with respect to a family of groups $\Cal H$ then each freely indecomposable free 
factor of $G$ must be included in $\Cal H$. In particular, a freely indecomposable 
Artin group is strongly relatively hyperbolic in only the most trivial of senses, namely with
respect to itself. 

In this paper we shall focus on \emph{weak} relative hyperbolicity and, 
following \cite{KS}, we shall consider hyperbolicity of an Artin group relative to its
finite type parabolic subgroups in this sense.  
By considering this question in connection with the geometry of
the associated Deligne complex (see Section \ref{Sect:Artin}) we are able to extend the
results of \cite{KS} considerably while giving a more unified approach to the problem.
Thus we have:

\begin{thm}\label{Thm:main} 
An Artin group is weakly hyperbolic relative  to its 
finite type standard parabolic subgroups if and only if its Deligne complex 
is a Gromov hyperbolic space.
\end{thm}

In Section \ref{Sect:2dim} we recall the definition of a \emph{2-dimensional} Artin group.
These include the groups considered in \cite{KS}. It is an interesting fact that the Deligne 
complex for a 2-dimensional Artin group is Gromov hyperbolic if and only if the Davis complex
associated with the underlying Coxeter group is. We therefore obtain the following 
(see Proposition \ref{Prop:2dim}).

\begin{thm}\label{Thm:2dim} 
A 2-dimensional Artin group is weakly hyperbolic relative to its 
finite type standard parabolic subgroups if and only if the corresponding Coxeter group
is a hyperbolic group.
\end{thm}

In the absence of a counter-example, it is reasonable to conjecture that the statement of 
Theorem \ref{Thm:2dim} holds for an arbitrary Artin group.

\begin{conj} \label{Conj:main}
An Artin group is weakly hyperbolic relative to its 
finite type standard parabolic subgroups if and only if the corresponding Coxeter group
is a hyperbolic group.
\end{conj}

In view of Theorem \ref{Thm:main} the conjecture 
is equivalent to stating that, for any Artin group,
the Deligne complex is Gromov hyperbolic if and only if the Davis complex associated 
to the corresponding Coxeter group is Gromov hyperbolic.  In one direction, this conjecture 
is easy to prove.  In Section \ref{Sect:2dim}, we give a simple proof that hyperbolicity of 
the Deligne complex implies hyperbolicity of the Davis complex. On the other hand, we are 
currently a long way from proving the converse.  The global geometry of the Deligne complex 
is not well understood in general.  It is not known, for example, whether every Deligne 
complex supports an equivariant metric of nonpositive curvature (a CAT(0) metric).  
Even in cases where such a metric known, the question of the existence of flat planes 
(and hence of Gromov hyperbolicity) is still rather delicate.
\medskip

Apart from the 2-dimensional case, the main situation where the Deligne complex 
is known to admit a CAT(0) metric is in the case of an Artin group of FC type. In this case the 
preferred metric is cubical. In Section \ref{Sect:FCtype}, we describe a systematic 
method of deforming the metric in each cube of the Deligne complex so as to obtain a
piecewise hyperbolic metric.  We than  obtain a  sufficient  condition  for this deformed 
metric to be CAT(-1). 
The condition is simply that the defining graph (as defined in Section \ref{Sect:Artin})
for the FC type Artin group has \emph{no empty squares}, meaning that
any circuit of length four has at least one diagonal pair of vertices spanning an edge.  
(We remark that ``no empty squares" in the defining graph does not imply ``no empty squares" 
in the link of every simplex of the Deligne complex, hence this deformation must be done 
carefully.) As a consequence we have the following.

\begin{thm}\label{Thm:FCtype} 
An FC type Artin group whose defining graph has no empty squares
is weakly hyperbolic relative to its finite type standard parabolic subgroups.
\end{thm}

Note that while all the Artin groups covered by Theorem \ref{Thm:FCtype} are necessarily 
associated with hyperbolic Coxeter groups, the ``no empty squares" condition does
not quite capture every hyperbolic Coxeter group of FC type. 
However, in this case, we are actually proving somewhat
more than is required for Theorem \ref{Thm:FCtype}, namely we show that the Deligne 
complex is CAT(-1) (which is a priori stronger than Gromov hyperbolic), 
and this for a rather specific choice of metric.

\medskip

We conclude this introduction with a remark on the proof of Theorem \ref{Thm:main}. 
The key observation is that the action of the Artin group on its Deligne complex is
particularly well-adapted to studying the group in relation to its finite type parabolic
subgroups simply because these subgroups are precisely the isotropy subgroups of the action.
One can therefore relate the geometry of the coned-off Cayley graph  to 
that of the Deligne complex (via a quasi-isometry) in order to prove  Theorem \ref{Thm:main}. 
As it happens, the most natural way of  expressing this argument is in the form of 
a more general ``relative" Milnor-Svarc Lemma. Recall that the usual Milnor-Svarc Lemma
states that if a finitely generated group $G$ acts properly discontinuously, cocompactly and
isometrically on a length space $X$, then $G$ and $X$ are quasi-isometric spaces. 
In particular, $G$ is a hyperbolic group if and only if $X$ is Gromov hyperbolic.
The relative version states that 
if a finitely generated group $G$  acts \emph{discontinuously} 
(i.e. with discrete orbits), cocompactly and isometrically on a  length
space $X$ then $X$ is quasi-isometric to the coned-off Cayley graph $\Ga_{S,\Cal H}$ for $G$, 
where $S$ is any finite generating set and $\Cal H$ denotes the collection of maximal 
isotropy subgroups. We defer the proof (and a careful statement) of this result to 
Section \ref{Sect:RMS}, Theorem \ref{Thm:RMS}, even though we shall use it almost immediately 
in Section \ref{Sect:Artin}, below.


\section{Relative hyperbolicity and Artin groups}\label{Sect:Artin}


Let $\Delta$ denote a simplicial graph with vertex set $V(\Delta)$ and 
edge set $E(\Delta)\subset V(\Delta)\times V(\Delta)$.
Suppose also that every edge $e=\{ s,t\}\in E(\Delta)$ 
carries a label $m_e=m_{st}\in \N_{\geq 2}$. 
We define the \emph{Artin group} $G(\Delta)$
associated to the (labelled) \emph{defining graph} $\Delta$
to be the group given by the presentation
\footnote{Our notion of defining graph differs from the frequently used
``Coxeter graph'' where, by contrast, the absence
of an edge between $s$ and $t$  indicates a commuting
relation ($m_{st}=2$) and the label $m_{st}=\infty$
is used to designate the absence of a relation between $s$ and $t$.
In our convention the label $\infty$ is never used.} 
\[
G(\Delta)=\<\ V(\Delta)\  \mid\  
\ubrace{ststs\cdots}{m_{st}}=\ubrace{tstst\cdots}{m_{st}}
\ \text{ for all } \{ s,t\}\in E(\Delta)\ \>\,.
\]
Adding the relations $s^2=1$ for each $s\in V(\Delta)$ yields a presentation 
of the associated \emph{Coxeter group} $W(\Delta)$ of type $\Delta$. We denote 
$\rho_\Delta:G(\Delta)\to W(\Delta)$ the canonical quotient map obtained by
this addition of relations. An Artin group is said to be of \emph{finite type}
(sometimes written \emph{spherical type}) if the associated Coxeter group is 
finite, and of \emph{infinite type} otherwise. 
By a \emph{standard parabolic subgroup} of $G(\Delta)$, or $W(\Delta)$, we mean 
any subgroup generated by a (possibly empty) subset of the standard generating set 
$V(\Delta)$. More generally, any subgroup which is conjugate to a standard 
parabolic subgroup (of $G(\Delta)$ or $W(\Delta)$) 
shall be referred to as a \emph{parabolic subgroup}.

Probably the most important tool currently used in the study of \emph{infinite} type 
Artin groups is the Deligne complex (see \cite{CD}, etc..). We described this 
complex in detail.

Consider $G=G(\Delta)$ for a fixed defining graph $\Delta$. For each subset
of the generating set, $R\subset V(\Delta)$, we shall write $\Delta_R$ for the 
full labelled subgraph of $\Delta$ spanned by $R$. (Here we attach a meaning to the
\emph{empty} defining graph $\Delta_\emptyset$ by setting 
$G(\Delta_\emptyset)=1$ and $W(\Delta_\emptyset)=1$).
The inclusion of $\Delta_R$ in $\Delta$ induces a
homomorphism $\phi_R : G(\Delta_R)\to G$ 
with image the standard parabolic subgroup $\<R\>$ generated by $R$. 
The construction of the Deligne complex is based on the rather nontrivial fact, due 
to H. van der Lek \cite{vdLek}, that, for every defining graph $\Delta$
and every $R\subset V(\Delta)$, the homomorphism $\phi_R$ is an isomorphism 
onto its image. Thus each standard parabolic subgroup of an Artin group is itself
canonically isomorphic to an Artin group. The corresponding 
statement for Coxeter groups is also true, and well-known (see \cite{Bou}).

Define
\[
\Cal V_f\ =\ \{ R\subset V(\Delta)\,:\ W(\Delta_R) \text{ finite } \}\,.
\]
We view $\Cal V_f$ as a partially ordered set under inclusion of sets, and 
define $K$ to be the geometric realisation of the derived complex of $\Cal V_f$. 
Thus there is a simplex $\sigma\in K$ of dimension $n\geq 0$ for every  
chain $R_0\subset R_1\subset\dots\subset R_n$ 
of $n+1$ distinct elements in $\Cal V_f$. We denote $\min(\sigma)=R_0$, 
the minimal vertex of $\sigma$.

Note that, for $\emptyset\subseteq R\subset T\subseteq V(\Delta)$, the inclusion 
$\Delta_R\subset\Delta_T$ induces a homomorphism 
$\phi_{R,T} : G(\Delta_R)\to G(\Delta_T)$. It follows that setting 
$G(\sigma)=G(\Delta_{\min(\sigma)})$ and 
$\phi_{\sigma,\tau}=\phi_{\min(\sigma),\min(\tau)}$, for all $\tau\subset\sigma\in K$,
defines a \emph{simple complex of groups} structure 
\[
(K,\{ G(\cdot)\},\{\phi_{\cdot,\cdot}\})\,
\]
in the sense of \cite{BH}. It is easily seen that the (orbifold) \emph{fundamental group}
of this complex of groups is isomorphic to the Artin group $G$ (via the homomorphisms
$\phi_R : G(\Delta_R)\to G$), and the result of van der Lek cited above ensures that the
complex of groups is \emph{developable}. It follows that the Artin group $G$ acts,
with quotient $K$ and isotropy subgroups the finite type parabolics
$\{\, G(\sigma)\, :\,\sigma\in K\,\}$, on a simply connected simplicial complex, 
the \emph{universal cover} of $(K,\{ G(\cdot)\},\{\phi_{\cdot,\cdot}\})$, which
we shall denote $\D$ and refer to as the \emph{Deligne complex} associated to $G(\Delta)$. 

We note that, replacing the collection of groups $\{\, G(\sigma)\, :\,\sigma\in K\,\}$
with the corresponding Coxeter groups $\{\, W(\sigma)\, :\,\sigma\in K\,\}$ 
leads in a similar way to the definition of a developable complex of groups 
whose fundamental group
is, this time, the Coxeter group $W=W(\Delta)$, and whose universal cover, which we shall
denote $\D_W$, is known as the \emph{Davis complex}. The Coxeter group 
acts on its Davis complex with finite vertex stabilizers (in fact properly discontinuously
and cocompactly). This is quite different from the situation of the Deligne complex
where the Artin group acts with every nontrivial vertex stabilizer an infinite group.
In particular, the Deligne complex is not even a locally compact space (while the Davis
complex clearly is). Nevertheless, a lot of important information is carried by the 
action of the Artin group on its Deligne complex, as can be seen from \cite{CD,Crisp} etc.
 
We suppose now that the complex $K$ is endowed with a piecewise Euclidean 
or piecewise hyperbolic metric. Since
$K$ is finite, this induces a complete $G$-equivariant length metric on the Deligne
complex (c.f. \cite{BH}). There are  two very natural choices for such a metric 
(namely the Moussong metric and the cubical metric) described in \cite{CD}.
These specific metrics are particularly useful when they can be shown to be 
nonpositively curved, or CAT(0), as demonstrated in \cite{CD}, where we refer the reader
for further details. However, in what follows, the actual choice of piecewise Euclidean 
or hyperbolic metric is more or less irrelevant since for any such metric $d$, the Deligne 
complex $(\D,d)$ will be ($G$-equivariantly) quasi-isometric to the 1-skeleton 
of $\D$ equipped with the unit-length edge metric. The statement ``the Deligne complex is 
Gromov hyperbolic'' shall henceforth be interpreted to mean ``with respect to any equivariant 
piecewise Euclidean or piecewise hyperbolic metric''.

We are now able to prove Theorem \ref{Thm:main} which we restate below:

\begin{thm}\label{main} 
An Artin group is weakly hyperbolic relative to its 
finite type standard parabolic subgroups if and only if its Deligne complex 
is a Gromov hyperbolic space.
\end{thm}

\begin{proof}
It is easily seen that the action of an Artin group on its Deligne complex (equipped with a 
piecewise Euclidean length metric) is discontinuous and co-compact. 
(As observed later, in Remark \ref{Ex:CofG},
this is actually true in the case of an arbitrary finite developable complex of groups).
Also, the (maximal) isotropy subgroups of this action are, by construction,
just the (maximal) finite type parabolic subgroups of the Artin group.
The result now follows immediately from the relative version of the Milnor-Svarc Lemma
(Theorem \ref{Thm:RMS}) proved in Section \ref{Sect:RMS}.
\end{proof}


\section{Two-dimensional Artin groups}\label{Sect:2dim}


We say that the Artin group $G=G(\Delta)$ is \emph{2-dimensional} if $\Delta$ has at
least one edge ($G$ is not free) and  every
triangle in $\Delta$ has edge labels $m,n,p$ satisfying $1/m+1/n+1/p\leq 1$, 
equivalently if every rank 3 parabolic subgroup is of infinite type. 
The terminology is justified by the fact that an Artin group is 2-dimensional 
in this sense if and only if it has cohomological dimension 2. Each 2-dimensional 
Artin group is also known to have geometric dimension 2. Moreover, the Deligne 
complex is 2-dimensional and is CAT(0) when equipped with the Moussong metric. 
These statements were all established in the paper of the first author and M.~Davis \cite{CD}.

In \cite{KS}, I. Kapovich and P. Schupp showed that an
Artin group with all indices $m_{ij}\geq 7$ is weakly hyperbolic relative 
 to its maximal finite type parabolic subgroups. 
It is clear that the groups treated by Kapovich and Schupp 
are all examples of 2-dimensional Artin groups. 
The following statement generalises their result.

\begin{prop}\label{Prop:2dim} 
Let $G(\Delta)$ be a 2-dimensional Artin group. Then the following are equivalent
\begin{itemize}
\item[(1)] the Davis complex $\D_W$ is Gromov hyperbolic;
\item[(2)] the Coxeter group $W(\Delta)$ is a hyperbolic group;
\item[(3)] $\Delta$ contains no triangle having edge labels $m,n,p$ with 
$1/m+1/n+1/p =1$ and no square with all edge labels equal to $2$. 
\item[(4)] the Deligne complex $\D$ (equipped with the Moussong metric) 
is Gromov hyperbolic;
\item[(5)] $G(\Delta)$ is weakly hyperbolic relative to 
the collection of finite type standard parabolic subgroups, namely the set
 $\{ G(e)\,:\, e\in E(\Delta)\,\}$.
\end{itemize}
\end{prop}

\begin{proof} The equivalence of (1) and (2) follows by an application of the
Milnor-Svarc Lemma and the equivalence of these and (3) follows immediately 
from Moussong's conditions \cite{Mou} (see discussion below). 
The equivalence of (1-3) and (4) is proved in \cite{Crisp} Lemma 5 
(using the Flat Plane Theorem and the relationship between the Deligne complex and the 
Davis complex). Finally, the equivalence of (4) and (5) is a consequence of the above
Theorem \ref{main}.
\end{proof}

We observe that there is quite generally a close relationship between the Deligne 
complex $\D$ associated to an Artin group $G=G(\Delta)$ and the Davis complex $\D_W$ 
for the corresponding Coxeter group $W=W(\Delta)$. 
The canonical projection induces a simplicial map $\D\to\D_W$ which is surjective,
and the Tits section (a setwise section to the canonical projection) induces 
an inclusion $\D_W\hookrightarrow \D$. We suppose that $\D_W$ is equipped 
with an equivariant metric. Since the Coxeter group $W$ acts properly co-compactly
and isometrically on $\D_W$ the space and the group are quasi-isometric
(by the Milnor-Svarc Lemma). Thus $W$ is a hyperbolic group precisely when $\D_W$ is
a Gromov hyperbolic space. 

By the work of Moussong \cite{Mou}, the Davis complex 
equipped with the Moussong metric $(\D_W, d_M)$ is, in all cases, 
a CAT(0) space and is Gromov hyperbolic if and only if it contains no isometrically embedded
flat plane.  
Moreover, Moussong shows that $\D_W$ contains an embedded flat only when $W$ contains an 
obvious $\Z \times \Z$ subgroup, namely when $W$ contains a standard parabolic subgroup 
which is either Euclidean or is the direct product of two infinite standard parabolics.

In summary, an arbitrary Coxeter group $W$ is hyperbolic group if and only if
both of the following conditions hold:
\begin{itemize}
\item[(M1)] no standard parabolic subgroup of $W$ is a Euclidean reflection group; and,
\item[(M2)] no standard parabolic subgroup of $W$ is the direct product of two 
infinite parabolics.
\end{itemize}

An analogous metric, also known as the Moussong metric, can be defined on 
the Deligne complex $\D$ such that the inclusion and projection maps  $\D_W \to \D \to \D_W$ 
are isometries on each simplex.  It is conjectured in [7] that the Moussong metric on $\D$ is 
always CAT(0).  This conjecture is only known to hold in the 
2-dimensional case and in some limited higher dimensional cases containing 
four-strand braid groups as parabolics (see \cite{Charn}).  
We can use this metric, however, to prove the easy direction of Conjecture \ref{Conj:main}.

\begin{prop}  
With respect to the Moussong metric,  the inclusion of $\D_W \hookrightarrow \D$ induced 
by the Tits section is an isometric embedding.  In particular, if $\D$ is  
Gromov hyperbolic, then so is $\D_W$.
\end{prop}

\begin{proof}
Since the natural projection $\D \to \D_W$ and the inclusion $\D_W \to \D$  take simplices 
isometrically onto simplices, they map geodesics to piecewise geodesics.  
It follows that both maps are distance non-increasing.  
But the composite map $\D_W \to \D  \to  \D_W$  is an isometry 
(in fact it is the identity map), so the first of these maps must also be distance 
non-decreasing, i.e., it is an isometric embedding.
\end{proof}


\section{FC type Artin groups}\label{Sect:FCtype}


The only other case in which the Deligne complex has been shown to have a CAT(0) metric 
is the case of  an Artin group of FC type.  A defining graph $\Delta$, and its associated 
Artin group $G(\Delta)$, is said to be of \emph{FC type} if it satisfies the following 
condition:
\begin{quote}
$T\subset V(\Delta)$ lies in $\Cal V_f$ if and only if $T$ spans a complete subgraph in $\Delta$.
\end{quote}
(Recall that a graph is complete if any two vertices are connected by an edge.)
There is a natural cubical structure on the Deligne complex of any Artin group 
which we will describe below.  The FC condition precisely guarantees that all 
links in this cubical structure are flag complexes (hence 
the terminology ``FC") and so the induced geodesic metric on $\D$ is CAT(0). 

For FC groups, there is a corresponding CAT(0) cubical structure on the 
Davis complex and the embedding $\D_W \hookrightarrow \D$ is isometric.  
It is thus clear that Moussong's conditions for hyperbolicity of the 
Davis complex are necessary conditions for hyperbolicity of the Deligne complex.
 We will prove a partial converse of this statement.   
We say that a graph $\Delta$ has \emph{no empty squares} if for any circuit of length four, 
at least one diagonal pair of vertices spans an edge.  We will 
show that if the defining graph for an Artin group is FC type and has no empty squares, 
then the Deligne complex $\D$ supports a hyperbolic metric. Note that these hypotheses 
are strictly stronger than Moussong's.  The FC condition rules out Euclidean parabolics 
in $W$, and the second condition rules out products of infinite parabolics (since by the FC 
condition, any infinite parabolic contains a pair of generators not connected by an edge in 
$\Delta$).  On the other hand, the graph below with all edges labelled 5 gives a defining 
graph satisfying both of Moussong's conditions (M1-M2), but neither of the two conditions above.

\begin{figure}[ht]
\begin{center}
\includegraphics[height=3cm]{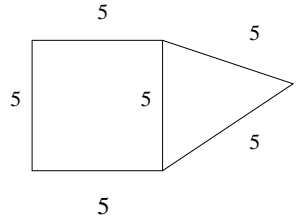}
\end{center}
\caption{Defining graph $\Gamma$} 
\label{Fig1}
\end{figure}

In the case of a right-angled Artin group, that is, one in which all edges in the 
defining graph are labeled $2$, the no empty squares condition exactly agrees with Moussong's 
conditions, thus Conjecture \ref{Conj:main} holds for these groups.

\medskip

The idea behind the proof of hyperbolicity is to slightly deform the CAT(0) cubical metric 
to make it negatively curved.  The no empty squares condition guarantees that some 
(but not all) of the links in this metric are ``extra large",  that is, the lengths 
of closed geodesics are bounded away from $2\pi$,  hence the metric at those vertices 
can be deformed slightly without destroying the CAT(0) condition.  The trick is to replace 
the Euclidean cubes with hyperbolic cubes which contain enough right angles to preserve the 
metric on those links which are not extra large.

\medskip

Define a hyperbolic metric on an $n$-cube as follows. Identify $n$-dimensional hyperbolic 
space with the hyperboloid 
\[
\mathbb H^n=\{ (x_1, \dots, x_n,x_{n+1}) \mid x_1^2 + \dots + x_{n-1}^2 - x_{n}^2 = -1, x_n>0\}
\]
in $\R^{n,1}$.  Let $H \cong (\Z/2)^n$ be the reflection group on $\mathbb H^n$ generated by 
reflections across the coordinate hyperplanes $x_i=0, i=1, \dots,n$. For any $\epsilon>0$, 
let $\xe$ be the (unique) point in the positive orthant of $\mathbb H^n$ at distance $\epsilon$ 
from every coordinate hyperplane. Let $Y^n_\epsilon$ be the convex hull in $\mathbb H^n$ of 
the $H$-orbit of $\xe$.  Then $Y^n_\epsilon$ is a regular hyperbolic $n$-cube of side 
length $2\epsilon$. Clearly, any $k$-dimensional face of $Y^n_\epsilon$ is isometric 
to $Y^k_\epsilon$.

Now let $\Ce$ be the intersection of $Y_\epsilon$ with the positive orthant.  
Then $\Ce$ is again a hyperbolic $n$-cube, but it is not regular.  
The vertex $x_0=(0, \dots, 0, 1)$ has all codimension one faces meeting at right 
angles whereas the vertex $\xe$ has all codimension one faces meeting at some 
angle $\theta<\frac{\pi}{2}$ that depends on $\epsilon$.  As $\epsilon$ goes to 
zero, the metric on $Y_\epsilon$ approaches the Euclidean metric, so we can 
take $\theta$ arbitrarily close to $\frac{\pi}{2}$  by choosing $\epsilon$ sufficiently small.  

Define the \emph{type} of a vertex $v$ in $\Ce$ to be the number of coordinate hyperplanes 
containing $v$.  Thus, $x_0$ has type $n$ while $\xe$ has type $0$.  A face $F$ of $\Ce$ 
contains a unique vertex of minimal type $k$, and a unique vertex of maximal type $l$.  
We say the face $F$ has \emph{type} $(k,l)$. By symmetry, it is easy to see that all 
faces of a given type $(k,l)$ in $\Ce$ are isometric. 
The next lemma shows that the metric on a face of type $(k,l)$ is independent of $n$.

\begin{lemma}
Suppose $F \subseteq \Ce$ is of type $(k,l)$.  Then $F$ is isometric to some (hence any) 
face of type $(k,l)$ in $C^l_\epsilon$. In particular, a face of type $(0,l)$ is isometric 
to $C^l_\epsilon$. 
\end{lemma}

\begin{proof}  
First consider a face $F$ of type $(0, l)$.  Any such face is the intersection 
of a $l$-dimensional face $\hat F$ of $Y^n_\epsilon$ with the positive orthant.  
Since $\hat F$ is isometric to $Y^l_\epsilon$, we conclude that $F$ is isometric 
to $C^l_\epsilon$.   More generally, any face $F$ of type $(k,l)$ is contained in
a face of type $(0,l)$, so $F$ is isometric to a face of type $(k,l)$ in $C^l_\epsilon$.
\end{proof}

Next, we analyze links of vertices in the cube $\Ce$.  By construction, the link 
of $x_0$ is a spherical simplex with all edge lengths $\pi/2$ (called an ``all-right 
spherical simplex"), while the link of $\xe$ is a regular spherical simplex with all 
edge lengths $\theta < \frac{\pi}{2}$.  

Suppose $v$ is a vertex of type $k$, $0<k<n$. Every codimension 1 face of $\Ce$ 
containing $v$ either  lies in a coordinate hyperplane, or in a codimension 1 face 
of $Y^n_\epsilon$.  Since a codimension 1 face of $Y^n_\epsilon$ is preserved by 
reflection across every coordinate  hyperplane it intersects, it follows that these 
two types of faces intersect orthogonally.  Thus,  the link of $v$ decomposes as an 
orthogonal join,
\[
link(v,\Ce)=link_\uparrow(v,\Ce) \ast link_\downarrow(v,\Ce)
\]
where $Link_\uparrow(v,\Ce) $ is the link of $v$ in the face $F_0$ spanned by $v$ 
and $x_0$, while $Link_\downarrow(v,\Ce)$ is the link of $v$ in the face $F_\epsilon$ 
spanned by $v$ and $\xe$.  The face $F_0$ is of type $(k,n)$ while $F_\epsilon$ is of 
type $(0,k)$.  By the lemma above, there is an isometry of $F_\epsilon$ with $C^k_\epsilon$ 
taking $v$ to the basepoint $x_0$ (in $C^k_\epsilon$). Thus, the downward link is an 
all-right simplex.

\medskip

We now define a piecewise hyperbolic metric on the Deligne complex $\D$.  
For this we will use a slightly different description of $\D$.  
Let $G=G(\Delta)$ and for $T \in \Cal V_f$, write $G_T=G(\Delta_T)$.  
A fundamental domain for the action of $G$ on $\D$ is the complex $K$ 
defined in the previous section.  If we think of the vertices of $K$ as 
finite type parabolic subgroups $G_T$ rather than as sets of generators $T\in \Cal V_f$, 
then the vertices in $\D$ correspond to cosets $gG_T$, $T \in \Cal V_f$ and simplices to 
totally ordered flags of cosets.

This gives rise to a cubical structure on $\D$ whose vertices are the cosets  $gG_T$ and 
whose cubes correspond to ``intervals",
\[
[gG_T, gG_{R}] = \textrm{span of the vertices $gG_{T'}$ 
with $gG_T \subseteq gG_{T' }\subseteq gG_{R}$}.
\]
Define an equivariant, piecewise hyperbolic metric $d_\epsilon$ on $\D$ by assigning each 
cube $[gG_T, gG_R]$ with $|T|=t$, $|R|=r$ the metric of a face of type $(t,r)$ 
in $C^r_\epsilon$.  In particular,  the cubes $[G_\emptyset, G_R]$ are isometric 
to $C^r_\epsilon$ with the vertex $G_\emptyset$ identified with $x_\epsilon$ 
and $G_R$ identified with  $x_0$.

\begin{thm}\label{FC}
Suppose  an Artin group $G(\Delta)$ is FC type, and its defining graph has no empty squares.
Then the metric $d_\epsilon$ on the Deligne complex $\D$ is CAT(-1)  for $\epsilon$ 
sufficiently small.  In particular, $\D$ is Gromov hyperbolic.
\end{thm}

\begin{proof} 
First consider the simplicail complex $L$ whose vertex set is $V(\Delta)$ and whose 
simplices $\sigma_T$ are spanned by the subsets $T \in \Cal V_f$.  (This is known as 
the ``link" of the Coxeter group $W(\Delta)$.) The hypotheses of the theorem precisely 
guarantee that this simplicial complex satisfies  Seibenman's ``no triangles, no squares" 
condition.  Moreover, the FC condition implies that the link of any simplex in $L$ is 
isomorphic (as a simplicial complex) to a full subcomplex of $L$.  Hence these links 
also satisfy the ``no triangles, no squares" condition.  In \cite{Gro1}, Gromov shows 
that if each simplex in such a complex is an all-right spherical simplex, then the 
resulting geodesic metric is ``extra large". That is, $L$ is CAT(1) and closed geodesics 
in $L$, and in all links in $L$, have lengths bounded away from $2\pi$.  As shown by 
Moussong in \cite{Mou}, any sufficiently small deformation of the metric  on each simplex 
preserves this property, hence the resulting metric remains CAT(1). 

To prove the theorem, we must show that the link of every vertex $v$ in $\D$ is CAT(1) 
with respect to the metric $d_\epsilon$.  By equivariance of the metric, it suffices to 
consider vertices $v=G_T$ lying in the fundamental domain $K$.  By the discussion above,
the link of $v$ in $\D$ decomposes as an orthogonal join
\[
link(v,\D)=link_\uparrow(v,\D) \ast link_\downarrow(v,\D),
\]
so it suffices to show that the upward and downward links are each CAT(1). 

The upward link consists of 
\[
link_\uparrow(v,\D)= \bigcup_{T \subset R} link(v, [G_T,G_R]).
\]
It has a $k$-simplex for each spherical parabolic $G_R \supset G_T$ with $|T |-|R|=k-1$.  
Thus, as an abstract simplicial complex, it can be identified with the complex $L$ defined 
above when $T=\emptyset$, and with  $link(\sigma_T, L)$ otherwise. The metric 
on $link_\uparrow(v,\D)$ induced by $d_\epsilon$ is a peicewise spherical metric 
with all edge lengths arbitrarily close to $\frac{\pi}{2}$ for sufficiently 
small $\epsilon$. It follows from the discussion above that this metric is CAT(1). 

 The downward link, $link_\downarrow(v,\D)$ is composed of the link of $G_T$ in the 
cubes $[gG_\emptyset,G_T]$, for $g \in G_T$.  It  is an all-right piecewise spherical 
complex which is isomorphic to the link of the maximal vertex $G_T$ in the Deligne 
complex for the finite type Artin group $G_T$.  
This was shown to be a flag complex in \cite{CD}, Lemma 4.3.2, hence it is CAT(1). 
 \end{proof}

\begin{cor}\label{FCrelhyp}
An Artin group of FC type whose defining graph has no empty squares is
 hyperbolic relative to its finite type standard parabolic subgroups.
\end{cor}


\section{A relative version of the Milnor-Svarc Lemma}\label{Sect:RMS}


Let $X$ be a metric space and $G$ a group which acts on $X$ by isometries.
We say that the action is \emph{co-compact} if there exists a compact set
$K\subset X$ such that $X=\bigcup\limits_{g\in G} gK$, and \emph{discontinuous} 
if every orbit is a discrete subspace of $X$, equivalently if 
for all $x\in X$ there exists an $\epsilon_x>0$ 
such that $d(x,y)>\epsilon_x$ for all $y\in G(x)\setminus\{ x\}$. 
A subgroup $H<G$ is said to be an 
\emph{isotropy subgroup} of $G$ if its fixed set in $X$ is non-empty.

\begin{thm}\label{Thm:RMS}
Let $G$ be a finitely generated group and suppose that $G$ admits 
a discontinuous, co-compact, isometric action on a length space $X$.
Let $\Cal H$ denote a collection of subgroups of $G$ consisting of exactly
one representative of each conjugacy class of maximal isotropy
subgroups for the action of $G$ on $X$.
Then $\Cal H$ is finite and, for any finite generating set $S$ of $G$, the
coned-off Cayley graph $\Gamma_{S,\Cal H}(G)$ is quasi-isometric to $X$.
In particular, if $X$ is a Gromov hyperbolic space then the group
$G$ is weakly hyperbolic relative to the collection $\Cal H$ of
maximal isotropy subgroups.
\end{thm}
 
\begin{proof}
We first show that the number of conjugacy classes of maximal isotropy
subgroups is finite. Clearly, distinct maximal isotropy subgroups have disjoint 
fixed sets in $X$. Also, since $G$ acts co-compactly, every maximal 
isotropy subgroup is conjugate to one which fixes a point inside a certain
compact region $K$ such that $X=\bigcup\limits_{g\in G} gK$.
By way of contradiction, we now suppose that there exists 
an infinite sequence $H_1,H_2,\dots$ of pairwise distinct maximal
isotropy subgroups such that each $H_i$ fixes a point $x_i\in K$ (the points $x_i$ 
being necessarily pairwise distinct). By compactness of $K$ we may pass to an infinite 
subsequence for which the sequence $(x_i)_{i\in\N}$ converges to a point $x\in K$.
Moreover, at most one of the $H_i$ may fix $x$, so we may as well suppose that 
none of them fix $x$. For each $i$ we may therefore choose an element 
$h_i\in H_i$ such that $h_i(x)\neq x$. 
Since $d(h_i(x),x)\leq 2d(x_i,x)$ it follows that the sequence
$(h_i(x))_{i\in\N}$ also converges to $x$, contradicting the assumption 
that $G$ acts discontinuously.
 
Let $H_1,H_2,\dots ,H_n$ denote the finitely many maximal isotropy subgroups whose
fixed sets intersect $K$ nontrivially. Set $Q=\{ g\in G : gK\cap K\neq\emptyset \text{ but }
g\notin H_r \text{ for all } r=1,\dots ,n\}$, and choose a subset $\what Q\subset Q$ which
contains exactly one representative for each coset $gH_r$, for 
$g\in Q$ and $r\in\{1,\dots ,n\}$. We note that $G$ is generated by the set $\what Q$ together
with the subgroups $H_1,\dots, H_n$. 

We now use a compactness argument to show that $\what Q$ is finite. 
By way of contradiction we suppose that there exists an infinite sequence 
$(g_i)_{i\in\N}$ of pairwise distinct elements of $\what Q$. For each $i\in\N$ we may find
$x_i,y_i\in K$ such that $g_i(x_i)=y_i\in gK\cap K$. Since $K\times K$ is compact we may pass
to an infinite subsequence of $(g_i)_{i\in\N}$ for which the sequence of pairs $(x_i,y_i)$
converges to a point $(x,y)\in K\times K$. Since 
\[
d(g_i(x),y)\leq d(g_i(x),g_i(x_i)) + d(g_i(x_i),y) = d(x,x_i)+d(y_i,y)
\]
it follows that the sequence $(g_i(x))_{i\in\N}$ converges to $y$.
Since the action of $G$ is discontinuous, this implies that 
the sequence is eventually constant: there exists $N$ such that $g_i(x)=y$ for all $i>N$.
But then, for any $N<i<j$, the element $g_j\inv g_i$ fixes $x\in K$ 
and hence lies in some maximal isotropy subgroup $H_r$. 
That is to say that $g_i$ and $g_j$ are
different representatives in $\what Q$ for the same coset of $H_r$, 
contradicting the choice of $\what Q$.

Since $G$ is finitely generated we may 
extract a finite generating set from any given generating set for the group.
It follows that we may extend $\what Q$ to a finite generating set 
$S$ of $G$ in such a way that 
$\what Q\subset S\subset \what Q\cup H_1\cup\dots\cup H_n$.
Let $\Ga=\Ga_{S, \Cal H}(G)$ denote the coned-off Cayley graph for
$G$ with respect to the generating set $S$ and a finite set $\Cal H$ of 
isotropy subgroups as prescribed in the statement of the Theorem. We may as well suppose 
that $\Cal H$ is chosen by selecting from amongst the subgroups $H_r$, $r=1,\dots, n$, 
one from each conjugacy class. Moreover, the exact choice of representatives 
for $\Cal H$ is not really important, as the structure of the coned-off Cayley graph $\Gamma$ 
depends only on the set of conjugacy classes of subgroups involved. 
Note also that, up to quasi-isometry, $\Ga$ depends only on the set $\what Q$ (rather 
than the choice of $S$) and the subgroups $H_r$. More generally, the coned-off Cayley graph
$\Ga=\Ga_{S, \Cal H}(G)$ is independent, up to quasi-isometry, of the choice of generating set
$S$ (regardless of whether or not it contains $\what Q$) as long as this set is finite.
    
Let $v_0$ denote the base vertex of $\Ga$ and write $\Ga_0$ for
the $G$-orbit of $v_0$ with metric induced from $\Ga$. Then the inclusion $\Ga_0\to \Ga$ 
is a quasi-isometry. We shall show that $\Ga_0$ (and therefore $\Ga$) is 
quasi-isometric to $X$.

Choose a point $x_0\in K$. 
This choice determines a $G$-equivariant map $f:\Ga_0\to X$ by sending 
$gv_0$ to $gx_0$ for all $g\in G$. 
It is easily seen that, for $p,q\in\Ga_0$, 
\[
d_X(f(p),f(q)) \leq R\, d_\Ga(p,q)\,,
\]
where 
$R$ denotes the maximum value in the finite set 
\[
\{\, d_X(x_0,s(x_0))\ :\ s\in\what Q\,\}
\cup\{\, 2d_X(x_0,\Fix(H_r))\ :\ r=1,\dots,n\,\}\,.
\]
In order to prove the reverse inequality (i.e. to bound $d_\Ga(p,q)$ above
by a linear function of $d_X(f(p),f(q))$) we need to establish the following fact:
\smallskip

\noindent 
\emph{There exists a constant $\epsilon >0$ such that, for all $g\in G$, either
$gK\cap K\neq\emptyset$ or $d_X(gK,K)>\epsilon$ (where here we understand the Hausdorff 
distance).} 
\smallskip

We use, once again, a compactness argument to prove this statement. 
If the statement is not true then we may find a sequence $(g_i)_{i\in\N}$
of distinct group elements such that 
$0<d_X(g_{i+1}K,K)<\frac{1}{2}d_X(g_iK,K)$ for all $i\in\N$.  
Choosing, for each $i$, a pair $(x_i,y_i)\in K\times K$ such that 
$d_X(g_iK,K)<d_X(g_i(x_i),y_i)<2d_X(g_iK,K)$, and passing to an infinite
subsequence for which the sequence of pairs converges 
to a pair $(x,y)\in K\times K$, we observe that 
the sequence $(g_i(x))_{i\in\N}$ converges to $y$. This contradicts
the assumption that the action of $G$ is discontinuous unless the sequence is eventually
stationary, that is, unless $g_i(x)=y$ for some $i$. But this is impossible since
$g_iK\cap K=\emptyset$ for all $i\in\N$.
\smallskip
  
We shall now give an upper bound for $d_\Ga(p,q)$. Since $X$ is a path metric space we
may find a path $\ga$ from $f(p)$ to $f(q)$ in $X$ whose length approximates 
the distance between these points to within $\epsilon$: 
$\ell(\ga)\leq d_X(f(p),f(q))+\epsilon$.
Choose $m\in\N$ such that $(m-1)\epsilon <\ell(\ga)\leq m\epsilon$, and
let $f(p)=y_0,y_1,\dots,y_m=f(q)$ denote equally spaced points along the path $\ga$.
In particular $d(y_{i-1},y_i)\leq \epsilon$, for all $i=1,\dots,m$. 
Choosing $K_0, K_1, K_2, \dots, K_m$ to be translates of the compact $K$ such that
$y_i\in K_i$ for all $i$, we observe, by the claim just proven, that
$K_{i-1}\cap K_i\neq\emptyset$ for all $i=1,\dots,m$.
By construction, whenever $gK\cap K\neq\emptyset$
we may express $g$ in the form $sh$ for $s\in \what Q$ and $h\in H_r$ (for some $r$).
Thus the sequence $K_0, K_1, K_2, \dots, K_m$ gives rise to a path of length
at most $2m$ joining $p$ to $q$ in $\Ga$. Thus $d_\Ga(p,q)\leq 2m$.
On the other hand $(m-1)\epsilon < \ell(\ga)\leq d_X(f(p),f(q))+\epsilon$.
Combining these inequalities results in
\[
d_\Ga(p,q)\ <\ \frac{2}{\epsilon}d_X(f(p),f(q))+4\,.
\]
This completes the proof that the map $f:\Ga_0\to X$ is a quasi-isometric embedding.
Clearly, since the compact $K$ is bounded, any point in $X$ is a bounded distance
from a point in the orbit of $x_0$, and so the map $f$ is in fact a quasi-isometry. 
\end{proof}

\subsection{Complexes of groups}
A rather general construction is to describe a group $G$ as the fundamental group of a
complex of groups (see \cite{BH}). If $G$ is the fundamental group of a finite 
complex $(Y,\Cal G)$ of groups which is developable and whose universal 
cover is Gromov hyperbolic, then $G$ is weakly hyperbolic relative to $\Cal G$. 
As discussed in the preceding section, each Artin group 
is the fundamental group of a finite complex of groups, 
leading to the results stated there.

\subsection{Mapping class groups}
A further example is that of the mapping class group $\Mod(S)$ of a closed 
orientable surface $S$ of higher genus. Recall that the complex of curves
$\Cal C(S)$ associated to the closed surface $S$ is defined to be the simplicial 
complex whose vertices are the nontrivial isotopy classes of simple closed curves 
and whose simplices are spanned by sets of vertices which are simultaneously represented 
by mutually disjoint (non-parallel) simple closed curves.
The group $\Mod(S)$ acts naturally on this complex (with rather large vertex stabilizers).
In \cite{MM}, Masur and Minsky showed that the complex of curves $\Cal C(S)$
is a Gromov hyperbolic space and used this to prove that the mapping class group 
$\Mod(S)$ is weakly relatively hyperbolic with respect to 
subgroups $H_C :=\{ g\in \Mod(S) : g[C]=[C]\}$ for a finite collection 
of simple closed curves $C$ in $S$. Their proof passes through a modified version 
of Teichmuller space (``electric space'') which they show to be quasi-isometric to 
$\Cal C(S)$. A proof of this result may also be obtained
by applying Theorem \ref{Thm:RMS} directly to the complex $\Cal C(S)$.
This avoids introducing the action on Teichmuller space while  
still invoking the hyperbolicity of the curve complex as proved in \cite{MM}.

\begin{ackn}
We would like to thank Mike Davis, Chris Hruska and 
Ilya Kapovich for a number of helpful remarks and observations in
connection with this work.
\end{ackn}

\end{document}